\documentclass[a4paper, conference]{IEEEtran}
\usepackage[T1]{fontenc}
\IEEEoverridecommandlockouts
\usepackage{cite}
\usepackage{amsmath,amssymb,amsfonts}
\usepackage{algorithmic}
\usepackage{graphicx}
\usepackage{textcomp}
\usepackage{xcolor}
\usepackage{mathtools, stackengine, mathrsfs, setspace, bm, multirow, tikz, pgfplots, adjustbox, calc, hyperref, url}
\usepackage{xparse}
\pgfplotsset{width=8.1cm, legend style={font=\small}}
\newsavebox{\fminipagebox}
\NewDocumentEnvironment{fminipage}{m O{\fboxsep}}
 {\par\kern#2\noindent\begin{lrbox}{\fminipagebox}
  \begin{minipage}{#1}\ignorespaces}
 {\end{minipage}\end{lrbox}%
  \makebox[#1]{%
    \kern\dimexpr-\fboxsep-\fboxrule\relax
    \fbox{\usebox{\fminipagebox}}%
    \kern\dimexpr-\fboxsep-\fboxrule\relax
  }\par\kern#2
 }
\begingroup\edef\x{\endgroup
  \mathchardef\mathdollar=\the\numexpr"7000+\the\mathdollar\relax
}\x
\DeclareMathAlphabet{\mathit}{T1}{cmr}{m}{it}

\makeatletter
\let\MYcaption\@makecaption
\makeatother
\usepackage{caption}
\DeclareCaptionLabelSeparator{colonquad}{:\quad}
\DeclareCaptionStyle{fig-top}%
	[justification=raggedright,indention=2cm, font=footnotesize, labelsep=colonquad]{}
\DeclareCaptionStyle{fig}%
	[justification=raggedright,indention=2cm, font=footnotesize, labelsep=colonquad]{}	
\makeatletter
\let\@makecaption\MYcaption
\makeatother

\pgfplotsset{
    every first x axis line/.style={},
    every first y axis line/.style={},
    every first z axis line/.style={},
    every second x axis line/.style={},
    every second y axis line/.style={},
    every second z axis line/.style={},
    first x axis line style/.style={/pgfplots/every first x axis line/.append style={#1}},
    first y axis line style/.style={/pgfplots/every first y axis line/.append style={#1}},
    first z axis line style/.style={/pgfplots/every first z axis line/.append style={#1}},
    second x axis line style/.style={/pgfplots/every second x axis line/.append style={#1}},
    second y axis line style/.style={/pgfplots/every second y axis line/.append style={#1}},
    second z axis line style/.style={/pgfplots/every second z axis line/.append style={#1}}
}

\makeatletter
\def\pgfplots@drawaxis@outerlines@separate@onorientedsurf#1#2{%
    \if2\csname pgfplots@#1axislinesnum\endcsname
    \else
    \scope[/pgfplots/every outer #1 axis line,
        #1discont,decoration={pre length=\csname #1disstart\endcsname, post length=\csname #1disend\endcsname}]
        \pgfplots@ifaxisline@B@onorientedsurf@should@be@drawn{0}{%
            \draw [/pgfplots/every first #1 axis line] decorate {
                \pgfextra
                \pgfplotspointonorientedsurfaceabsetupfor{#2}{#1}{\pgfplotspointonorientedsurfaceN}%
                \pgfplots@drawgridlines@onorientedsurf@fromto{\csname pgfplots@#2min\endcsname}%
                \endpgfextra 
                };
        }{}%
        \pgfplots@ifaxisline@B@onorientedsurf@should@be@drawn{1}{%
            \draw [/pgfplots/every second #1 axis line] decorate {
                \pgfextra
                \pgfplotspointonorientedsurfaceabsetupfor{#2}{#1}{\pgfplotspointonorientedsurfaceN}%
                \pgfplots@drawgridlines@onorientedsurf@fromto{\csname pgfplots@#2max\endcsname}%
                \endpgfextra 
                };
        }{}%
    \endscope
    \fi
}%
\makeatother

\def\BibTeX{{\rm B\kern-.05em{\sc i\kern-.025em b}\kern-.08em
    T\kern-.1667em\lower.7ex\hbox{E}\kern-.125emX}}

\begin{document}
\bstctlcite{IEEEexample:BSTcontrol}
\title{Quantification of the Impact of \textit{GHG} Emissions\\ on Unit Commitment in Microgrids
\thanks{This work was supported in part by the Research Council of Norway under the ``LUCS'' project, and by the German Federal Ministry for Economic Affairs and Energy under Grant 03EI6004B.}
}

\author{\IEEEauthorblockN{Ogun Yurdakul,
Fikret Sivrikaya and
Sahin Albayrak}
\IEEEauthorblockA{Department of Electrical Engineering and Computer Science\\
Technische Universit{\"a}t Berlin, Berlin, Germany}
Email: \{yurdakul, fikret.sivrikaya, sahin.albayrak\}@tu-berlin.de\\\\
\small{\textit{presented in IEEE PES T\&D-LA 2020}}
\vspace{-0.5cm}}

\maketitle

\begin{abstract}
The global climate change creates a dire need to mitigate greenhouse gas (\textit{GHG}) emissions from thermal generation resources (\textit{TGR}s). While microgrids are instrumental in enabling the deeper penetration of renewable resources, the short-term planning of microgrids needs to explicitly assess the full range of impact of \textit{GHG} emissions. To this end, we propose a novel unit commitment (\textit{UC}) approach, which enables the representation of \textit{GHG} emissions from \textit{TGR}s, the specification of \textit{GHG} emission constraints, and the ex-ante evaluation of carbon tax payment with all other costs and benefits. We quantify the relative merits of the proposed approach vis-à-vis the classical \textit{UC} approach via representative studies. The results indicate that the proposed \textit{UC} approach yields lower costs than does the classical \textit{UC} approach and achieves a greater reduction in costs as carbon tax rate increases. Further, increasing carbon tax rates can markedly disincentivize \textit{TGR} generation under the proposed approach. 
\end{abstract}
\section{Introduction}\label{1}
A microgrid is a cluster of loads, distributed generation resources (\textit{DGR}s), and electric storage resources (\textit{ESR}s) that operate in coordination to supply electricity in a reliable manner. Typically integrated to its host power system at the distribution level, a microgrid is perceived by its distribution system as a single entity responding to appropriate signals \cite{eucbib:101}. For all intents and purposes, a microgrid is a microcosm of a bulk power system that retains most of its innate operational characteristics.\par
The \textit{DGR}s in a microgrid can be broadly bifurcated into two categories: thermal generation resources (\textit{TGR}s) and variable energy resources (\textit{VER}s). \textit{TGR}s include microturbines, fuel cells, and reciprocating internal combustion engines with generators and are especially common in microgrids in rural areas, developing nations, and military premises \cite{eucbib:69}. \textit{TGR}s have controllable power output but can undergo only gradual temperature changes and hence are subject to minimum uptime, minimum downtime, and ramping constraints \cite{eucbib:63}.\par
\textit{VER}s, such as photovoltaic (\textit{PV}) panels and wind turbines, are characterized by a renewable fuel source that can be neither stored nor controlled. \textit{VER}s cannot be similarly situated to \textit{TGR}s, since \textit{VER} power outputs are highly time-varying, intermittent, and uncertain. Further, microgrids with integrated \textit{ESR}s, such as batteries, ultracapacitors, and flywheels, possess various capabilities, including the hallmark capability to store electric energy for later use.\par
Similar to bulk power systems, the short-term planning of a microgrid can be determined via unit commitment (\textit{UC}) and economic dispatch (\textit{ED}) decisions \cite{eucbib:79}. The classical unit commitment (\textit{CUC}) approach seeks minimum cost strategies to determine the start-up and shut-down of \textit{TGR}s based on expected load, equipment limitations, and operational policies \cite{eucbib:63}. The equipment limitations of \textit{TGR}s and the inter-temporal constraints of microgrid physical asset operations render \textit{UC} a time-coupled problem, and necessitate that the \textit{UC} decisions be taken typically one-hour to one-week ahead of operations based on the uncertain data/information available at the time of decision. The injections levels of \textit{TGR}s are subsequently determined by the solution of the \textit{ED} problem after most of the uncertainty unravels.\par
Microgrids lend themselves as conducive environments to enabling the deeper penetration of \textit{VER}s without unduly exacerbating the stress on transmission systems. As such, microgrids, aided by \textit{UC} approaches that evaluate the full breadth of impact of greenhouse gas (\textit{GHG}) emissions, can play a pivotal rule in combating global climate change.\par
The thorough assessment of \textit{GHG} emissions in short-term microgrid operations hinges on practical \textit{UC} approaches that expressly include \textit{GHG} emission models. Such approaches further need to have the capability to stipulate explicit constraints on the amount of \textit{GHG} emissions over a study period.\par
Another major requirement, brought on especially with the advent of carbon pricing schemes, is the analysis of the monetary impacts of \textit{GHG} emissions. Carbon tax sets a specific price on the amount of emitted carbon dioxide equivalent ($CO_{2}e$) to internalize the negative externalities of \textit{TGR} generation. The prevalence of carbon pricing has been soaring in recent years, with the price reaching $\mathit{\$}139/tCO_{2}e$ in Sweden \cite{eucbib:62}. \par
\begin{table*}
\begin{minipage}[c][11.0cm][t]{\linewidth}
\begin{fminipage}{\textwidth}[1ex]
\footnotesize
\renewcommand{\arraystretch}{1.2}
\begin{tabular}{l l  l l }
\multicolumn{2}{l}{\textbf{Nomenclature}}  & & \\
$\mathscr{H}$ & set of simulation time periods &\hspace{0.1cm} $[p^{\mathsf{w}}_{\sigma_s}]^m$ & minimum charging power of \textit{ESR} $\sigma_s$ (\textit{kW}) \\
$\mathscr{G}$ & set of distributed generation resources (\textit{DGR}s)& \hspace{0.1cm} $[p^{\mathsf{w}}_{\sigma_s}]^M$ & maximum charging power of \textit{ESR} $\sigma_s$ (\textit{kW})\\
$\mathscr{G}_{\text{\textit{VER}}}$ & set of variable energy resources (\textit{VER}s)& \hspace{0.1cm} $[E_{\sigma_s}]^{m}$ & minimum energy storage limit of \textit{ESR} $\sigma_s$ (\textit{kWh})\\
$\mathscr{G}_{\text{\textit{TGR}}}$ & set of thermal generation resources (\textit{TGR}s)&\hspace{0.1cm} $[E_{\sigma_s}]^{M}$ & maximum energy storage limit of \textit{ESR} $\sigma_s$ (\textit{kWh})\\
$\mathscr{S}$ & set of electric storage resources (\textit{ESR}s)&\hspace{0.1cm} $\eta^{\mathsf{i}}_{\sigma_s} \in (0,1]$ & discharging efficiency of \textit{ESR} $\sigma_s$ \\ 
${h}$ & index of an hourly time period &\hspace{0.1cm} $\eta^{\mathsf{w}}_{\sigma_s} \in (0,1]$  & charging efficiency of \textit{ESR} $\sigma_s$\\
${\gamma_g}$ & distributed generation resource $g$  &\hspace{0.1cm} $p^{\mathsf{w}}_{\delta}[h]$& total microgrid load in hour $h$ (\textit{kW}) \\
$\sigma_s$ & electric storage resource $s$ &\hspace{0.1cm} $\lambda[h]$ & the price at which the microgrid purchases (resp. sells) \\
$[p^{\mathsf{i}}_{\gamma_g} ]^{m}$ & minimum power output of \textit{DGR} $\gamma_g$ (\textit{kW}) & \hspace{0.1cm} & energy from (resp. to) the distribution company in\\
$[p^{\mathsf{i}}_{\gamma_g} ]^{M}$ & maximum power output of \textit{DGR} $\gamma_g$ (\textit{kW})&\hspace{0.1cm} &hour $h$ ($\mathit{\$}/kWh$) \\
$[T^{\uparrow}_{\gamma_g} ]^{m}$ & minimum uptime of \textit{TGR} $\gamma_g$ (\textit{hrs})&\hspace{0.1cm}  $\big[R[h]\big]^{m}$ & spinning reserve requirement for hour $h$ (\textit{kW})\\
$[T^{\downarrow}_{\gamma_g} ]^{m}$ & minimum downtime of \textit{TGR} $\gamma_g$ (\textit{hrs})&\hspace{0.1cm}  $p_{\gamma_g}^{\mathsf{i}}[h]$ & power generation of \textit{DGR} $\gamma_g$ in hour $h$ (\textit{kW})\\
$\overline{\overline{c}}_{\gamma_g}$, $\overline{c}_{\gamma_g}$, $c_{\gamma_g}$ & quadratic ($\mathit{\$}/kW^2h$), linear ($\mathit{\$}/kWh$), and fixed& \hspace{0.1cm}  $u^{\mathsf{i}}_{\gamma_g}[h]$ & commitment status of \textit{TGR} $\gamma_g$ in hour $h$\\
&  ($\mathit{\$}/h$) fuel cost parameter of \textit{TGR} $\gamma_g$& \hspace{0.1cm} $r_{\gamma_g}^{\mathsf{i}}[h]$ & spinning reserve of \textit{TGR} $\gamma_g$ in hour $h$\\
$\mu_{\gamma_g}$ & start-up cost of \textit{TGR} $\gamma_g$ ($\mathit{\$}$)& \hspace{0.1cm} $u^{\mathsf{i}}_{\sigma_s}[h]$ & injection status of $\sigma_s$ in hour $h$\\
$\overline{\overline{k}}_{\gamma_g}$, $\overline{k}_{\gamma_g}$, $k_{\gamma_g}$ & quadratic ($kgCO_2e/kW^2h$), linear ($kgCO_2e/kWh$), & \hspace{0.1cm} $u^{\mathsf{w}}_{\sigma_s}[h]$ & withdrawal status of $\sigma_s$ in hour $h$\\
& and fixed ($kgCO_2e/h$) \textit{GHG} emission parameter of  & \hspace{0.1cm} $p^{\mathsf{i}}_{\sigma_s}[h]$ & power injection of $\sigma_s$ in hour $h$ (\textit{kW})\\
& \textit{TGR} $\gamma_g$ & \hspace{0.1cm} $p^{\mathsf{w}}_{\sigma_s}[h]$ & power withdrawal of $\sigma_s$ in hour $h$ (\textit{kW})\\
$[\kappa]^M$ & maximum \textit{GHG} emission limit for the study period  & \hspace{0.1cm}  $p^{\mathsf{n}}_{\sigma_s}[h]$ & net power injection of $\sigma_s$ in hour $h$ (\textit{kW})\\
&($kgCO_2e$) & \hspace{0.1cm}  $E_{\sigma_s}[h]$ & energy stored in $\sigma_s$ in hour $h$ (\textit{kWh})\\
$\psi$ & carbon tax rate ($\mathit{\$/kgCO_2e}$) & \hspace{0.1cm} $p_{\varphi}^{\mathsf{n}}[h]$ & net power injection of the distribution system to the\\
$\varkappa_{\gamma_g}$ & \textit{GHG} emissions of \textit{TGR} $\gamma_g$ over the study period & &microgrid in hour $h$ (\textit{kW})\\
& ($kgCO_2e$)& \hspace{0.1cm} $\xi^{\dagger}_{\gamma_g} $ & fuel cost of \textit{TGR} $\gamma_g$ over the study period ($\mathit{\$}$)\\
$\kappa $ & total \textit{GHG} emissions from all microgrid \textit{TGR}s over the & \hspace{0.1cm} $\xi^{\ddagger}_{\gamma_g} $ & total start-up cost of \textit{TGR} $\gamma_g$ over the study period ($\mathit{\$}$)\\
&study period ($kgCO_2e$)& \hspace{0.1cm} $\xi_{\varphi} $ & total net cost for the exchange of power with the\\
$[p^{\mathsf{i}}_{\sigma_s}]^m$ & minimum discharging power of \textit{ESR} $\sigma_s$ (\textit{kW}) & &distribution company ($\mathit{\$}$)\\
$[p^{\mathsf{i}}_{\sigma_s}]^M$ & maximum discharging power of \textit{ESR}  $\sigma_s$ (\textit{kW})& \hspace{0.1cm} $\xi_{\varkappa}  $ & carbon tax payment ($\mathit{\$}$) \\
\end{tabular}
\end{fminipage}
\end{minipage}
\end{table*}
The \textit{CUC} approach does not consider \textit{GHG} emissions and so does not take into account the carbon tax payment due to the \textit{GHG} emissions from \textit{TGR}s at the time of decision. As such, the carbon tax payment, for which the microgrid is liable, is evaluated \textit{ex-post} under the \textit{CUC} approach, thereby potentially bringing about dire economic implications. The \textit{UC} approaches for microgrids must conduct an \textit{ex-ante} evaluation of the carbon tax payment that will have been incurred, as per the effective carbon pricing schemes. The ex-ante evaluation of carbon tax payment permits a more thorough quantification of the benefits of taking \textit{UC} decisions that favor the greater utilization of \textit{VER}s jointly with \textit{ESR}s in lieu of \textit{TGR}s.\par
While economic mechanisms can serve as prime movers for major change, a key issue that needs to be investigated is whether carbon tax rates can effectively deter the use of \textit{TGR}s and incentivize the further utilization of \textit{VER}s in conjunction with \textit{ESR}s. As such, the analytical study of the influence of carbon pricing schemes on \textit{UC} decisions can provide useful guidance for policy makers addressing global warming.\par
\subsection{Related Work}
There is a growing body of literature on \textit{UC} approaches for microgrids. In \cite{eucbib:79}, the authors propose a \textit{UC} approach for microgrids with integrated \textit{TGR}s, \textit{VER}s, and \textit{ESR}s, yet the proposed approach does not consider the \textit{GHG} emissions from integrated \textit{TGR}s or evaluate the monetary impacts of \textit{GHG} emissions. The frameworks presented in \cite{eucbib:33, eucbib:19} consider \textit{ESR}s, \textit{TGR}s, and \textit{VER}s in the \textit{UC} of a microgrid; nonetheless, they do not model the \textit{GHG} emissions from \textit{TGR}s or their economic implications.\par
In \cite{eucbib:95}, the \textit{UC} decisions for microgrids with integrated \textit{TGR}s and \textit{VER}s have been studied, where the \textit{GHG} emissions from \textit{TGR}s are modeled. However, \cite{eucbib:95} does not include an explicit constraint on the amount of \textit{GHG} emissions, study carbon tax rates, or include \textit{ESR}s---which are key to enable the greater utilization of \textit{VER}s and so to mitigate \textit{GHG} emissions. While the approach presented in \cite{eucbib:96} models the \textit{GHG} emissions from \textit{TGR}s of a microgrid, it does not impose any constraints on the amount of \textit{GHG} emissions or evaluate the impact of carbon tax rates. In \cite{eucbib:104}, the \textit{UC} of a microgrid is studied, where the \textit{GHG} emissions from \textit{TGR}s and carbon emission costs are modeled. However, the authors of \cite{eucbib:104} do not study the influence of carbon tax rate on \textit{UC} decisions or model an explicit constraint for \textit{GHG} emissions from \textit{TGR}s in the \textit{UC} problem formulation.
\subsection{Contributions of the paper}
The general contributions and novel aspects of this paper are as follows:
\begin{enumerate}
\item We propose a novel \textit{UC} approach that comprehensively evaluates the full range of impact of \textit{GHG} emissions. We refer to this approach as \textit{environmental unit commitment} or \textit{EUC}. To the best of our knowledge, this is the first approach that simultaneously models \textit{GHG} emissions, allows the stipulation of a constraint on \textit{GHG} emission amount, and conducts an ex-ante evaluation of the carbon tax payment. We conduct representative studies and demonstrate the effectiveness of the proposed \textit{EUC} approach on real-world data\footnote{The source code is available at: https://github.com/oyurdakul/euc}. 
\item We examine the sensitivity of the \textit{UC} decisions to carbon tax rate and study the extent to which carbon tax can influence the operation of microgrid \textit{TGR}s, \textit{ESR}s, and power exchange with the distribution grid. Our findings regarding the influence of carbon tax rate on short-term microgrid operations could prove useful for policy makers in judiciously determining carbon tax rate. 
\end{enumerate}
This paper contains four additional sections. In Section \ref{2}, we develop models for the microgrid physical, economic, and environmental aspects. In Section \ref{3}, we present the mathematical formulations of the \textit{EUC} and \textit{CUC} approaches. We illustrate the capabilities and effectiveness of the proposed \textit{EUC} approach in Section \ref{4} with representative studies and discuss the results. We summarize the paper and provide directions for future work in Section \ref{5}. 
\section{Microgrid Modeling}\label{2}
We devote this section to the delineation of the microgrid models. We discretize the time axis and adopt one hour as the smallest indecomposable unit of time. In line with \cite{eucbib:98}, we decompose the study period into $H$ non-overlapping hours and 
define the study period by the set $\mathscr{H} \coloneqq \{h \colon h=1,...,H\}$. 
\subsection{Physical Asset Models}
We consider a microgrid interfaced with the distribution system and denote by $\mathscr{G}$ the set of \textit{DGR}s in the microgrid. We define the subsets $\mathscr{G}_{\text{\textit{VER}}}$ and $\mathscr{G}_{\text{\textit{TGR}}}$ to denote the set of \textit{VER}s and \textit{TGR}s, respectively, and we write the relation $\mathscr{G}=\mathscr{G}_{\text{\textit{VER}}}\bigcup\mathscr{G}_{\text{\textit{TGR}}}$. We define by $p_{\gamma_g}^{\mathsf{i}}[h]$ the \textit{kW} power injection of $\gamma_g \in \mathscr{G}$ in hour $h$. \par
The binary variable $u^{\mathsf{i}}_{\gamma_g}[h]\in\{0,1\}$ denotes the commitment status of \textit{TGR} $\gamma_g \in \mathscr{G}_{\text{\textit{TGR}}}$ in hour $h$. $u^{\mathsf{i}}_{\gamma_g}[h]=1$ if $\gamma_g$ is up in hour $h$, and $0$ otherwise. We define by $r_{\gamma_g}^{\mathsf{i}}[h]$ the spinning reserve of \textit{TGR} $\gamma_g \in \mathscr{G}_{\text{\textit{TGR}}}$ in hour $h$. We denote by $p^{\mathsf{w}}_{\delta}[h]$ the total \textit{kW} load of the microgrid in hour $h$.\par
We consider that an \textit{ESR} may inject power, withdraw power, or remain idle in hour $h$. Let $u^{\mathsf{w}}_{\sigma_s}[h]=1$ if the \textit{ESR} $\sigma_s$ withdraws power in hour $h$, and $0$ otherwise. Similarly, let $u^{\mathsf{i}}_{\sigma_s}[h]=1$ if $\sigma_s$ injects power in hour $h$, and $0$ otherwise. We denote by $p^{\mathsf{i}}_{\sigma_s}[h]$ and $p^{\mathsf{w}}_{\sigma_s}[h]$ the power injection and power withdrawal of $\sigma_s$ in hour $h$, respectively. We define by $p^{\mathsf{n}}_{\sigma_s}[h] \coloneqq p^{\mathsf{i}}_{\sigma_s}[h]-p^{\mathsf{w}}_{\sigma_s}[h]$ the net power injection of $\sigma_s$ in hour $h$. Let $E_{\sigma_s}[h]$denote the energy stored in the $\sigma_s$ at the end of hour $h$, or equivalently, at the beginning of hour $h+1$.\par
We assume that the distribution company (\textit{DisCo}) is the sole owner and operator of the distribution system with which the microgrid is interfaced. Let $p_{\varphi}^{\mathsf{n}}[h]$ denote the \textit{net} power injection of the distribution system to the microgrid in hour $h$. We adopt the convention that if the distribution system injects (resp. withdraws) power to (resp. from) the microgrid in hour $h$, then $p_{\varphi}^{\mathsf{n}}[h]>0$ (resp. $p_{\varphi}^{\mathsf{n}}[h]<0$). If there is no exchange of power between the microgrid and the distribution system in hour $h$, then $p_{\varphi}^{\mathsf{n}}[h]=0$.
\subsection{Environmental Models}\label{2b}
We devote this subsection to the development of the models for \textit{GHG} emissions from the microgrid \textit{TGR}s. We explicitly represent the total $kgCO_2e$ \textit{GHG} emission of each \textit{TGR} $\gamma_g$ over the study period by the relation
\begin{equation}
\label{varkappa}
\varkappa_{\gamma_g}= \sum\limits_{ h \in \mathscr{H}} \Big[ \big(\overline{\overline{k}}_{\gamma_g}(p_{\gamma_g}^{\mathsf{i}}[h])^2+\overline{k}_{\gamma_g}(p_{\gamma_g}^{\mathsf{i}}[h])+k_{\gamma_g}\big)u_{\gamma_g}^{\mathsf{i}}[h] (1 \text{\textit{ hr}}) \Big], 
\end{equation}
based on the \textit{GHG} emission modeling in \cite{eucbib:95}. The terms $\overline{\overline{k}}_{\gamma_g}$, $\overline{k}_{\gamma_g}$, and $k_{\gamma_g}$ in \eqref{varkappa} denote the quadratic ($kgCO_2e/kW^2h$), linear ($kgCO_2e/kWh$), and fixed ($kgCO_2e/h$) \textit{GHG} emission parameter of the \textit{TGR} $\gamma_g$, respectively.\par
We further express the total $kgCO_2e$ \textit{GHG} emissions from all microgrid \textit{TGR}s over the study period by 
\begin{equation}
\label{kappa}
\kappa=\sum_{\gamma_g \in \mathscr{G}_{\text{\textit{TGR}}}} \varkappa_{\gamma_g}.
\end{equation}
The proposed \textit{EUC} approach allows the stipulatation of an upper limit on the total \textit{GHG} emissions from microgrid \textit{TGR}s over the study period. Such an upper limit ensures that, independent of the economic factors, the generation of microgrid \textit{TGR}s is explicitly constrained by the resulting \textit{GHG} emissions. To this end, we model by $[\kappa]^M$ the maximum \textit{GHG} emission limit for the study period.
\subsection{Economic Models}
In this subsection, we model the costs and benefits associated with the microgrid operation over the study period. We express the fuel cost of \textit{TGR} $\gamma_g$ over the study period by
\begin{equation}
\xi^{\dagger}_{\gamma_g}= \sum\limits_{h \in \mathscr{H}}\Big[ \big(\overline{\overline{c}}_{\gamma_g}(p_{\gamma_g}^{\mathsf{i}}[h])^2+	\overline{c}_{\gamma_g}(p_{\gamma_g}^{\mathsf{i}}[h])+c_{\gamma_g}\big)u_{\gamma_g}^{\mathsf{i}}[h]  (1 \text{\textit{ hr}})\Big], \label{fcost} \\
\end{equation}
based on \cite{eucbib:79}. We express the total start-up cost of \textit{TGR} $\gamma_g$ over the study period by  
\begin{equation}
\xi^{\ddagger}_{\gamma_g}= \sum\limits_{h \in \mathscr{H}}\Big[\mu_{\gamma_g}(1-u^{\mathsf{i}}_{\gamma_g}[h-1])u^{\mathsf{i}}_{\gamma_g}[h]\Big]. \label{sucost} 
\end{equation}\par
We consider that the \textit{DisCo} uses \textit{time-of-use rates} and utilizes \textit{net metering} as the billing mechanism. We express the total net cost (\textit{i.e.}, total cost minus total benefit) associated with the exchange of power with the \textit{DisCo} by
\begin{equation}
\label{net}
\xi_{\varphi}=\sum\limits_{h \in \mathscr{H}} \Big[\lambda[h]p^{\mathsf{n}}_{\varphi}[h] (1 \text{\textit{ hr}}) \Big]. 
\end{equation}\par
The \textit{CUC} approach evaluates only $\xi^{\dagger}_{\gamma_g}$ and $\xi^{\ddagger}_{\gamma_g}$  as the costs associated with the operation of \textit{TGR} $\gamma_g$. A key source of cost that is not captured by the \textit{CUC} approach is the carbon tax payment associated with the \textit{GHG} emissions from \textit{TGR}s over the study period. Our objective is to enable the ex-ante evaluation of carbon tax payment simultaneously with $\xi^{\dagger}_{\gamma_g}$ and $\xi^{\ddagger}_{\gamma_g}$. As such, we expressly model the carbon tax payment in the proposed \textit{EUC} approach.\par
We denote by $\psi$ the \textit{carbon tax rate}, which is the price for each unit of $kgCO_2e$ emitted, evaluated in $\mathit{\$}/kgCO_2e$. Utilizing $\psi$, as well as the total \textit{GHG} emissions from the microgrid \textit{TGR}s $\kappa$ modeled in Subsection \ref{2b}, we express the carbon tax payment for the \textit{GHG} emissions from microgrid \textit{TGR}s over the study period by the relation
\begin{equation}
\xi_{\varkappa} = \psi \kappa.
\end{equation}
\section{\textit{EUC} Problem Formulation}\label{3}
In this section, we present the \textit{EUC} problem formulation using the physical asset models, environmental models, and economic models developed in Section \ref{2}.
The general \textit{EUC} problem formulation is stated as:
\begin{equation}
\label{obj}
\begin{aligned}&\textit{EUC:} & 
\underset{\begin{subarray}{c}
{\scriptscriptstyle u^{\mathsf{i}}_{\gamma_g}[h], p^{\mathsf{i}}_{\gamma_g}[h], } \\
{\scriptscriptstyle u^{\mathsf{i}}_{\sigma_s}[h], p^{\mathsf{i}}_{\sigma_s}[h], } \\
{\scriptscriptstyle u^{\mathsf{w}}_{\sigma_s}[h], p^{\mathsf{w}}_{\sigma_s}[h], }\\
{\scriptscriptstyle p^{\mathsf{n}}_{\varphi}[h], r_{\gamma_g}^{\mathsf{i}}[h]}
\end{subarray}}{\text{minimize}} &\;\sum_{\gamma_g \in \mathscr{G}_{\text{\textit{TGR}}}} \Big[\xi^{\dagger}_{\gamma_g} + \xi^{\ddagger}_{\gamma_g}\Big] + \xi_{\varphi}+\xi_{\varkappa} ,  
\end{aligned}
\end{equation}
\hspace{18mm}subject to
{\allowdisplaybreaks
\begin{IEEEeqnarray}{l}
u^{\mathsf{i}}_{\gamma_g}[h] [p^{\mathsf{i}}_{\gamma_g}]^{m}   \leq   p^{\mathsf{i}}_{\gamma_g}[h]    \leq  u^{\mathsf{i}}_{\gamma_g}[h] [p^{\mathsf{i}}_{\gamma_g}]^{M}, \label{c1} \\
p^{\mathsf{i}}_{\gamma_g}[h] +r^{\mathsf{i}}_{\gamma_g}[h]    \leq  u^{\mathsf{i}}_{\gamma_g}[h] [p^{\mathsf{i}}_{\gamma_g}]^{M}, \label{c2} \\
u^{\mathsf{i}}_{\gamma_g}[h]  -  u^{\mathsf{i}}_{\gamma_g}[h-1]  \leq  u^{\mathsf{i}}_{\gamma_g}[\nu],\, \forall \nu \in \mathbb{N}\: \text{such that} \nonumber\\
\quad \quad \quad \quad h \leq \nu \leq \textrm{min}\{h-1+T^{\uparrow}_{\gamma_g},H\}, \label{c3}  \\
u^{\mathsf{i}}_{\gamma_g}[h-1] - u^{\mathsf{i}}_{\gamma_g}[h]  \leq 1- u^{\mathsf{i}}_{\gamma_g}[\nu],\, \forall \nu \in \mathbb{N}\: \text{such that}\; \nonumber\\
\quad \quad \quad \quad h \leq \nu  \leq \textrm{min}\{h-1+T^{\downarrow}_{\gamma_g},H\}, \label{c4}  \\
u^{\mathsf{i}}_{\sigma_s}[h]   +  u^{\mathsf{w}}_{\sigma_s}[h]   \leq 1, \label{c5}\\
u^{\mathsf{i}}_{\sigma_s}[h]  [p^{\mathsf{i}}_{\sigma_s}]^{m}  \leq  p^{\mathsf{i}}_{\sigma_s}[h]    \leq  u^{\mathsf{i}}_{\sigma_s}[h] [p^{\mathsf{i}}_{\sigma_s}]^{M}, \label{c6}\\
u^{\mathsf{w}}_{\sigma_s}[h] [p^{\mathsf{w}}_{\sigma_s}]^{m}  \leq  p^{\mathsf{w}}_{\sigma_s}[h]   \leq  u^{\mathsf{w}}_{\sigma_s}[h] [p^{\mathsf{w}}_{\sigma_s}]^{M}, \label{c7}\\
E_{\sigma_s}[h]  =  E_{\sigma_s}[h-1]+ \frac{1}{\eta^{\mathsf{w}}_{\sigma_s}} p^{\mathsf{w}}_{\sigma_s}[h] -\eta^{\mathsf{i}}_{\sigma_s} p^{\mathsf{i}}_{\sigma_s}[h], \label {c8}\\
{[E_{\sigma_s}]}^{m} \leq   E_{\sigma_s}[h] \leq  [E_{\sigma_s}]^{M},  \label{c9}\\
\sum_{\gamma_g \in \mathscr{G}} p_{\gamma_g}^{\mathsf{i}}[h] + \sum_{\sigma_s \in \mathscr{S}} p_{\sigma_s}^{\mathsf{n}}[h] + p_{\varphi}^{\mathsf{n}}[h] = p^{\mathsf{w}}_{\delta}[h], \label{c10} \\
\sum_{\gamma_g \in \mathscr{G}_{\text{\textit{TGR}}}} r^{\mathsf{i}}_{\gamma_g}[h] \geq \big[R[h]\big]^{m}, \label{c11}\\
\kappa\leq [\kappa]^M, \label{c12}
\end{IEEEeqnarray}}
where we take into account the constraints \eqref{c1}-\eqref{c11} $\forall h \in \mathscr{H}$, the constraints \eqref{c1}-\eqref{c4} $\forall \gamma_g \in \mathscr{G}_{\text{\textit{TGR}}}$, and the constraints \eqref{c5}-\eqref{c9} $\forall \sigma_s \in \mathscr{S}$.\par
The \textit{EUC} approach seeks to minimize the carbon tax payment, jointly with the fuel and start-up costs of \textit{TGR}s, and the net cost associated with the exchange of power with the \textit{DisCo}, as expressed by the objective function \eqref{obj}. The \textit{EUC} problem formulation explicitly considers the constraints on \textit{TGR} outputs by \eqref{c1} and \eqref{c2}. \textit{TGR} minimum uptime and minimum downtime constraints are expressed by \eqref{c3} and \eqref{c4}, respectively. We represent the constraint that \textit{ESR}s may not both inject and withdraw power at the same time by \eqref{c5}. The constraints on the power injection and withdrawal of \textit{ESR}s are expressed by \eqref{c6} and \eqref{c7}, respectively. We express the intertemporal operational constraint of \textit{ESR}s by \eqref{c8}. The constraints on the energy stored in \textit{ESR}s are expressed by \eqref{c9}. The power balance of the microgrid is ensured by \eqref{c10}. We state the constraints on spinning reserve requirements by \eqref{c11}. The constraint \eqref{c12} stipulates a limit on the allowable \textit{GHG} emissions over the study period, which may effectively restrict the generation of \textit{TGR}s.\par
To illustrate the quantitative improvements imparted by the proposed \textit{EUC} approach, we present the \textit{CUC} problem formulation, which does not consider the \textit{GHG} emissions at the time of decision. \textit{CUC} problem formulation is stated as:
\begin{equation}
\label{obj2}
\begin{aligned}&\textit{CUC:} & 
\underset{\begin{subarray}{c}
{\scriptscriptstyle u^{\mathsf{i}}_{\gamma_g}[h], p^{\mathsf{i}}_{\gamma_g}[h], } \\
{\scriptscriptstyle u^{\mathsf{i}}_{\sigma_s}[h], p^{\mathsf{i}}_{\sigma_s}[h], } \\
{\scriptscriptstyle u^{\mathsf{w}}_{\sigma_s}[h], p^{\mathsf{w}}_{\sigma_s}[h], }\\
{\scriptscriptstyle p^{\mathsf{n}}_{\varphi}[h],  r_{\gamma_g}^{\mathsf{i}}[h]}
\end{subarray}}{\text{minimize}} \;\sum_{\gamma_g \in \mathscr{G}_{\text{\textit{TGR}}}} \Big[\xi^{\dagger}_{\gamma_g} + \xi^{\ddagger}_{\gamma_g}\Big] + \xi_{\varphi},  
\end{aligned}
\end{equation}
\hspace{21mm}subject to \hspace{4mm} \eqref{c1}-\eqref{c11}.

\section{Case Study and Results}\label{4}
In this section, we illustrate the application and effectiveness of the proposed \textit{EUC} approach on representative studies. 
\subsection{Case Study Data}
We consider a microgrid connected to the low-voltage side of the distribution transformer to power residential loads. We consider that the microgrid includes a diesel generator denoted by $\gamma_1$, a \textit{PV} panel denoted by $\gamma_2$, and an \textit{ESR} denoted by $\sigma_1$. The peak load of the microgrid is $37$ \textit{kW}.\par
The $\gamma_1$ has a peak capacity of $50$ \textit{kW}. The $\gamma_1$ cost parameters are $\overline{\overline{c}}_{\gamma_1}=\mathit{\$}1.20(10^{-3})/kW^2h$, $\overline{c}_{\gamma_1}=\mathit{\$}0.208/kWh$, and $c_{\gamma_1}=\mathit{\$}3.2/h$. Further, the $\gamma_1$ \textit{GHG} emission parameters are $\overline{\overline{k}}_{\gamma_1}=3.03(10^{-3})\,kgCO_2e/kW^2h$, $\overline{k}_{\gamma_1}=0.53\,kgCO_2e/kWh$, and $k_{\gamma_1}=8.09\,kgCO_2e/h$. The data for $\gamma_1$ and $\sigma_1$ are extracted from \cite{eucbib:79, eucbib:83} and presented in Table \ref{tgr_esr_table}. The peak capacity of $\sigma_2$ is 17 \textit{kW}. The load and \textit{PV} generation data are extracted from \cite{eucbib:80} and contain measurements for an anonymous house in New York. Since the load and \textit{PV} data are collected in New York, to ensure consistency, we consider the time-of-use rates offered by Con Edison, \textit{viz.}: 21.97\textit{\textcent/kWh} from 8 a.m. to midnight and 1.55\textit{\textcent/kWh} from midnight to 8 a.m. \cite{eucbib:89}. We consider a $15\%$ minimum reserve requirement measured with respect to the microgrid peak load over the study period provided solely by \textit{TGR} $\gamma_1$. We explicitly stipulate a constraint on the allowable \textit{GHG} emissions over the study period and take $[\kappa]^{M}=220\,kgCO_2e$. We further consider that the carbon tax rate is $\psi=\mathit{\$}0.07/kgCO_2e$. 
\begin{table}[h]
\vspace{-0.3cm}
\footnotesize
\renewcommand{\arraystretch}{1.5}
\caption{Case study data for \textit{TGR} $\gamma_1$ and \textit{ESR} $\sigma_1$}
\label{tgr_esr_table}
\centering
\begin{tabular}{c | c | | c | c | | c | c }
\hline \hline
parameter & value & parameter & value & parameter & value \\
\hline \hline
$[p^{\mathsf{i}}_{\gamma_1}]^{m}$ & $5$ \textit{kW} & $[p^{\mathsf{i}}_{\gamma_1}]^{M}$ & $50$ \textit{kW} & $\mu_{\gamma_1}$ & $\mathit{\$}1$  \\
$[T^{\uparrow}_{\gamma_1}]^{m}$ & $2$ \textit{hrs} & $[T^{\downarrow}_{\gamma_1}]^{M}$ & $2$ \textit{hrs} & $u^{\mathsf{i}}_{\gamma_1}[0]$ & 0 \\
$[p^{\mathsf{i}}_{\sigma_1}]^{m}$ & 0 \textit{kW} & $[p^{\mathsf{w}}_{\sigma_1}]^{m}$ & $0$ \textit{kW} & $[E_{\sigma_1}]^{m}$ & $0$ \textit{kWh} \\
$[p^{\mathsf{i}}_{\sigma_1}]^{M}$ & $12$ \textit{kW} & $[p^{\mathsf{w}}_{\sigma_1}]^{M}$ & $12$ \textit{kW} & $[E_{\sigma_1}]^{M}$ & $30$ \textit{kWh}  \\
\hline \hline
\end{tabular}
\vspace{-0.3cm}
\end{table} 
\subsection{Load and \textit{PV} Forecasting }
The computation of the numerical solutions of the \textit{EUC} and \textit{CUC} problems requires the numerical representation of \textit{VER} generation and microgrid load over the study period. To this end, in this section, we utilize the methodology presented in \cite{eucbib:88} to forecast \textit{PV} generation and microgrid load over the study period. The utilized methodology leverages a sequence-to-sequence (\textit{S2S}) architecture that comprises two long short-term memory (\textit{LSTM}) networks, \textit{viz.}: encoder and decoder. \par
We construct one \textit{S2S} architecture for each of the two forecasting tasks. We provide each \textit{S2S} architecture with the measurements for the previous 24 hours as well as the hour of the day and the day of the week of the forecasted time periods. Further, each \textit{S2S} architecture generates forecasts for the subsequent 24 hours. It is worth emphasizing that, since the data in \cite{eucbib:80} were anonymized and the exact location of the house was not disclosed, we did not provide the \textit{S2S} architectures with relevant weather data, such as cloudiness index or temperature. Nevertheless, different studies can utilize other forecasting methodologies, probability distributions, or historical data in the implementation of the \textit{EUC} approach.\par
\begin{figure}[!h]
\vspace{-0.3cm}
\centering
\pgfplotsset{scaled x ticks=false}
\vspace{0cm}
\begin{tikzpicture}
\pgfplotsset{}
	\begin{axis}[ scaled ticks=false, tick label style={/pgf/number format/fixed},
width=0.5\textwidth,
height=0.32\textwidth,
xlabel={\textit{hour }$h$},
ylabel={\textit{power (kW)}},
xlabel style={at={(axis description cs:0.5, 0)},anchor=north},
ylabel style={at={(axis description cs:0.1,.5)},rotate=0,anchor=south},
xmin=1,
xmax=24,
xtick={1,6,12,18,24},
ymin=0,
ytick={0,5,10,15,20,25,30},
ymax=30,
legend style={font=\small, at={(0.48,-0.27)}, legend columns=2, anchor=north},
label style={font=\small},
tick label style={font=\small}]
	\addplot+[black!30!green, smooth, solid, semithick, mark=triangle] plot coordinates
		{ (1,0.36) (2,0.64) (3,0.59) (4,0.53) (5,0.82) (6,1.25) (7,1.99) (8,4.05) (9,6.88) (10,9.11) (11,11.40) (12,12.39) (13,13.76) (14,13.84) (15,12.57) (16,10.74) (17,8.29) (18,5.50) (19,2.14) (20,1.30) (21,1.24) (22,1.09) (23,1.02) (24,1.07) };
	\addlegendentry{\textit{forecasted} $\gamma_2$ \textit{generation}  \hspace{0.34cm}}
	\addplot+[black!70!green, densely dashed, semithick, mark=x] plot coordinates
		{ (1,0) (2,0) (3,0) (4,0) (5,0) (6,0.36) (7,0.91) (8,2.08) (9,3.72) (10,6.73) (11,9.12) (12,11.08) (13,13.76) (14,12.74) (15,11.83) (16,11.78) (17,6.58) (18,3.29) (19,2.62) (20,0.05) (21,0) (22,0) (23,0) (24,0) };
	\addlegendentry{\textit{actual} $\gamma_2$ \textit{generation}  \hspace{0.34cm}}
	\addplot+[white!10!orange, smooth, solid, semithick, mark=diamond] plot coordinates
		{ (1,11.86) (2,10.63) (3,9.43) (4,9.00) (5,8.40) (6,11.13) (7,14.25) (8,14.64) (9,14.10) (10,15.45) (11,9.47) (12,14.57) (13,12.85) (14,17.52) (15,19.25) (16,27.56) (17,29.75) (18,24.83) (19,26.00) (20,24.45) (21,21.93) (22,15.79) (23,11.81) (24,10.05) };
	\addlegendentry{\textit{forecasted microgrid load}   \hspace{0.34cm}}
	\addplot+[black!30!orange, smooth, densely dashed, semithick, mark=+] plot coordinates
		{ (1,8.67) (2,8.57) (3,7.70) (4,8.46) (5,8.37) (6,9.82) (7,9.59) (8,13.98) (9,13.71) (10,12.39) (11,12.61) (12,16.10) (13,17.54) (14,16.68) (15,24.44) (16,23.26) (17,26.73) (18,29.96) (19,26.99) (20,24.18) (21,19.59) (22,19.62) (23,14.37) (24,10.85) };
	\addlegendentry{\textit{actual microgrid load}  \hspace{0.34cm}}
	\end{axis}
\end{tikzpicture}
\vspace{-0.3cm}
\caption{Forecasted and actual $\gamma_2$ generation and microgrid load }
\vspace{-0.1cm}
\label{forresults}
\end{figure}
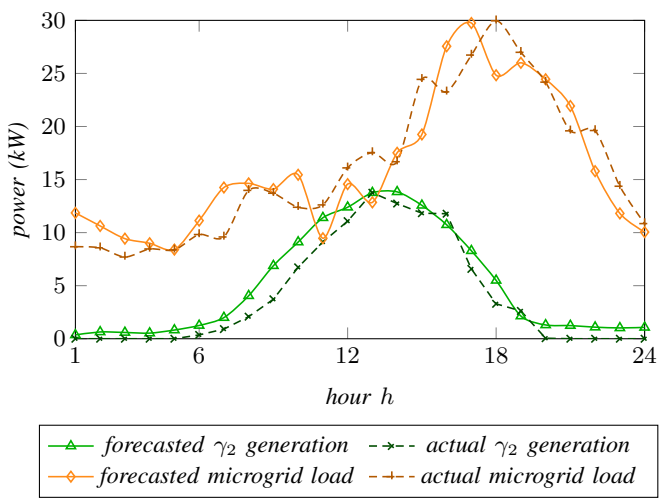
The dataset of each \textit{S2S} architecture contains measurements collected between May 1, 2019 and July 29, 2019 at one-hour resolution. Each dataset is split into training (60\%), validation (20\%), and test (20\%) sets, and we use the validation sets to tune the hyperparameters of the \textit{S2S} architectures. We pick the date of July 29, 2019 as the study period. The values for July 28, 2019 and July 29, 2019 are in the test set; therefore, the \textit{S2S} architectures have not been trained with the specific study period or the preceding day values. We utilize Tensorflow 
and Keras
to train, validate, and test the \textit{S2S} architectures on an NVIDIA Tesla P100 16 \textit{GB} \textit{GPU} with 800 \textit{GB} of \textit{RAM}.\par
The \textit{S2S} architecture to forecast \textit{PV} generation contains one \textit{LSTM} layer comprising 512 \textit{LSTM} blocks in both the encoder and decoder networks. We use the \textit{Adam} optimizer, and to prevent the networks from overfitting, utilize dropout with a probability of 0.4. The \textit{S2S} architecture to forecast \textit{PV} generation yields an \textit{RMSE} of 0.2876 on the test set.\par
The \textit{S2S} architecture to forecast microgrid load has one \textit{LSTM} layer comprising 512 \textit{LSTM} blocks in the encoder network and two \textit{LSTM} layers each comprising 512 \textit{LSTM} blocks in the decoder network. We pick \textit{Adam} as the optimizer and use dropout with a probability of 0.5. The \textit{S2S} architecture to forecast the microgrid load achieves an \textit{RMSE} of 0.5825 on the test set. Fig. \ref{forresults} depicts the \textit{S2S} architecture forecasts for the study period along with the corresponding actual measurements. The \textit{S2S} architecture forecasts for July 29, 2019 are utilized in the \textit{EUC} and \textit{CUC} solutions to represent $\gamma_2$ generation and microgrid load over the study period. 
\subsection{Unit Commitment Results}
\begin{figure}[!h]
\vspace{-0.3cm}
\centering
\begin{tikzpicture}
\pgfplotsset{}
\begin{axis}[
xlabel={\textit{hour }$h$},
ylabel={\textit{power (kW)}},
width=0.5\textwidth,
height=0.32\textwidth,
xlabel style={at={(axis description cs:0.5,0)},anchor=north},
ylabel style={at={(axis description cs:0.1,.5)},rotate=0, anchor=south},
xmin=1,
xmax=24,
xtick = {1,6,12,18,24},
ymin=-16,
ymax=16,
ytick = {-16,-12,-8,-4,0,4,8,12,16},
legend style={font=\small, at={(0.5,-0.27)}, legend columns=3, anchor=north},
label style={font=\small},
tick label style={font=\small}  ]
	\addplot+[white!10!olive, smooth, solid, semithick, mark=o] plot coordinates
		{ (1,5.0) (2,5.0) (3,5.0) (4,5.0) (5,5.0) (6,5.0) (7,5.0) (8,5.0) (9,6.724) (10,6.724) (11,6.724) (12,6.724) (13,6.724) (14,6.724) (15,6.724) (16,6.724) (17,6.724) (18,6.724) (19,6.724) (20,6.724) (21,6.724) (22,6.724) (23,6.724) (24,6.724) };
	\addlegendentry{\textit{EUC} $p_{\gamma_1}^{\mathsf{i}}[h]$\hspace{0.12cm}\,\hspace{0.12cm}}
	\addplot+[black!30!red, smooth, solid, semithick, mark=triangle] plot coordinates
		{ (1,-4.265) (2,-3.750) (3,-3.554) (4,-3.467) (5,-3.467) (6,-3.599) (7,-3.811) (8,-4.088) (9,2.970) (10,2.481) (11,2.075) (12,1.774) (13,1.568) (14,1.431) (15,1.341) (16,1.295) (17,1.292) (18,1.334) (19,1.419) (20,1.549) (21,1.744) (22,2.037) (23,2.493) (24,3.196) };
	\addlegendentry{\textit{EUC} $p_{\sigma_1}^{\mathsf{n}}[h]$\hspace{0.12cm}\,\hspace{0.12cm}}
	\addplot+[white!10!cyan, smooth, solid, semithick,mark=diamond] plot coordinates
		{ (1,10.765) (2,8.740) (3,7.394) (4,6.937) (5,6.047) (6,8.479) (7,11.071) (8,9.678) (9,-2.474) (10,-2.865) (11,-10.729) (12,-6.318) (13,-9.202) (14,-4.475) (15,-1.385) (16,8.801) (17,13.444) (18,11.272) (19,15.717) (20,14.877) (21,12.222) (22,5.939) (23,1.573) (24,-0.940) };
	\addlegendentry{\textit{EUC} $p_{\varphi}^{\mathsf{n}}[h]$\hspace{0.12cm}\,\hspace{0.12cm}}
	\addplot+[black!60!olive, smooth, semithick, mark=+, densely dashed] plot coordinates
{ (1,5.0) (2,5.0) (3,5.0) (4,5.0) (5,5.0) (6,5.0) (7,5.0) (8,5.0) (9,9.458) (10,9.458) (11,9.458) (12,9.458) (13,9.458) (14,9.458) (15,9.458) (16,9.458) (17,9.458) (18,9.458) (19,9.458) (20,9.458) (21,9.458) (22,9.458) (23,9.458) (24,9.458) };
	\addlegendentry{\textit{CUC} $p_{\gamma_1}^{\mathsf{i}}[h]$\hspace{0.12cm}\,\hspace{0.12cm}}
		\addplot+[black!60!red, smooth, solid, semithick, mark=x] plot coordinates
		{ (1,-4.265) (2,-3.750) (3,-3.554) (4,-3.467) (5,-3.467) (6,-3.599) (7,-3.811) (8,-4.088) (9,2.970) (10,2.481) (11,2.075) (12,1.774) (13,1.568) (14,1.431) (15,1.341) (16,1.295) (17,1.292) (18,1.334) (19,1.419) (20,1.549) (21,1.744) (22,2.037) (23,2.493) (24,3.196) };
	\addlegendentry{\textit{CUC} $p_{\sigma_1}^{\mathsf{n}}[h]$\hspace{0.12cm}\,\hspace{0.12cm}}
		\addplot+[black!30!blue, smooth, solid, semithick, mark=star, densely dashed] plot coordinates
		{ (1,10.765) (2,8.749) (3,7.394) (4,6.937) (5,6.047) (6,8.479) (7,11.071) (8,9.678) (9,-5.208) (10,-5.599) (11,-13.464) (12,-9.052) (13,-11.937) (14,-7.209) (15,-4.120) (16,6.067) (17,10.709) (18,8.538) (19,12.983) (20,12.142) (21,9.488) (22,3.205) (23,-1.162) (24,-3.674) };
	\addlegendentry{\textit{CUC} $p_{\varphi}^{\mathsf{n}}[h]$\hspace{0.12cm}\,\hspace{0.12cm}}
	\end{axis}
\end{tikzpicture}
\vspace{-0.4cm}
\caption{Optimal operations under the \textit{EUC} and \textit{CUC} approaches}
\vspace{-0.1cm}
\label{res}
\end{figure}
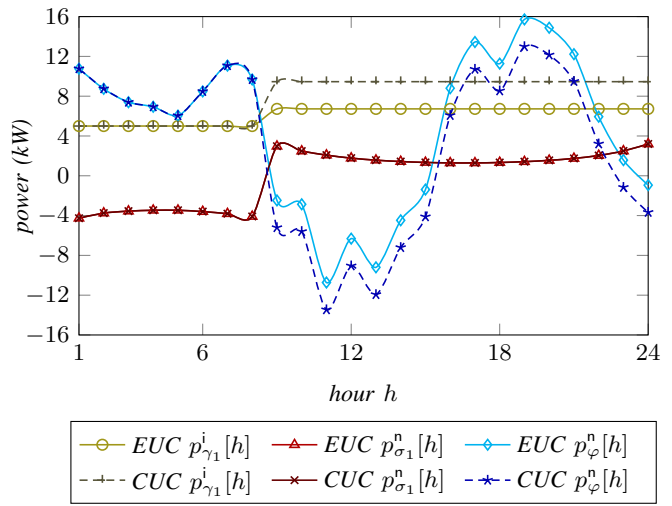

The \textit{EUC} problem formulation described by \eqref{obj}-\eqref{c12} is a mixed-integer-programming (\textit{MIP}) problem known to be \textit{NP}-hard. We solve the \textit{EUC} problem using Gurobi 8.1 
on a 2.6 \textit{GHz} Intel Core i7 \textit{CPU} with 16 \textit{GB} of \textit{RAM} for the study period.
Fig. \ref{res} presents the optimal injections and withdrawals under the \textit{EUC} and \textit{CUC} approaches. In both formulations, $\sigma_1$ tends to charge (resp. discharge) when the \textit{DisCo} electricity rates are low (resp. high), thereby exploiting intraday price variation and capitalizing on arbitrage opportunities.\par
The marked difference between the \textit{EUC} and \textit{CUC} approaches manifests itself in the optimal operations from hour 9 to hour 24. The optimal \textit{EUC} solution generates less energy from the \textit{TGR} compared to the optimal \textit{CUC} solution, since the \textit{EUC} approach is cognizant of the carbon tax payment while taking \textit{UC} decisions. On the flip side, the \textit{CUC} approach does not consider the carbon tax payment while taking \textit{UC} decisions and so needs to conduct an ex-post evaluation of the carbon tax payment for which the microgrid is liable. Owing to this ex-ante evaluation of carbon tax payment, the total costs for the study period are 16.9\textit{\textcent} lower under the \textit{EUC} approach than those under the \textit{CUC} approach. \par
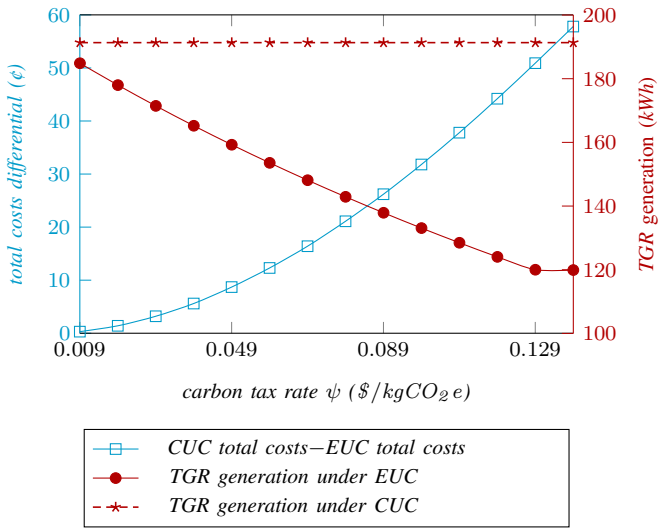
\begin{figure}[!h]
\vspace{0cm}
\centering
\begin{tikzpicture}

\pgfplotsset{
  xmin=0.009,
  xmax=0.139,
  xtick = {0.009,0.049,0.089,0.129},
  x tick label style={/pgf/number format/.cd,
            fixed,
            fixed zerofill,
            precision=3,
        /tikz/.cd},
        label style={font=\footnotesize},
  	tick label style={font=\footnotesize},
        scaled ticks=false, tick label style={/pgf/number format/fixed},}

	\begin{axis}[ 
	axis y line*=left,
	y axis line style=black!20!cyan,
	ytick style=black!20!cyan,
	label style={font=\footnotesize},
  tick label style={font=\footnotesize},
  width=0.445\textwidth,
height=0.32\textwidth,
  xlabel={ \textit{carbon tax rate} $\psi$ \textit{(}$\mathit{\$/kgCO_2e}$\textit{)}},
  ylabel={\textit{total costs differential (\textcent)}},
  y tick label style ={black!20!cyan},
  xlabel style={at={(axis description cs:0.5, 0)},anchor=north},
  ylabel style={at={(axis description cs:0.1,.5)},rotate=0,anchor=south, black!20!cyan},
  ymin=0,
  ymax=60,
  ytick={0,10,20,30,40,50,60},
    ]
	\addplot+[black!20!cyan, smooth, mark=square]  coordinates
   		{ (0.009,0.3) (0.019,1.4) (0.029,3.2) (0.039,5.6) (0.049,8.7) (0.059,12.3) (0.069,16.4) (0.079, 21.1) (0.089,26.2) (0.099,31.8) (0.109,37.8) (0.119,44.2) (0.129,50.9) (0.139,57.8) };
		\label{plot1}
	\end{axis}

\begin{axis}[
  axis y line*=right,
  axis x line=none,
    width=0.445\textwidth,
height=0.32\textwidth,
  y axis line style=black!30!red,
  ylabel style={at={(axis description cs:1.38,.5)},rotate=0,anchor=south, black!30!red},
  ymin=100, ymax=200,
  legend style={ at={(0.5,-0.3)}, font=\footnotesize, legend columns=1, anchor=north},
  ylabel=\textit{TGR} generation (\textit{kWh}), 
  y tick label style=black!30!red,
  ytick style=black!30!red,
  label style={font=\footnotesize},
  tick label style={font=\footnotesize},
]
\addplot+[black!20!cyan, smooth, mark=square]
coordinates
   		{ (0.009,0.3)};
\addlegendentry{\hspace{0.3cm}\textit{CUC total costs}$-$\textit{EUC total costs}\hspace{0.9cm}\quad \hspace{1.6cm}}
\addplot[smooth,mark=*,black!30!red]
   coordinates
   		{ (0.009,184.86) (0.019,177.99) (0.029,171.46) (0.039,165.22) (0.049,159.26) (0.059,153.57) (0.069,148.12) (0.079, 142.90) (0.089,137.89) (0.099,133.08) (0.109,128.47) (0.119,124.04) (0.129,120) (0.139,119.9) };
\addlegendentry{\hspace{0.3cm}\textit{\textit{TGR} generation under \textit{EUC}}\hspace{1.5cm}\quad \hspace{1.6cm}}
\addplot[densely dashed, semithick, mark=star,black!30!red]
  coordinates
   { (0.009,191.33) (0.019, 191.33) (0.029,191.33) (0.039,191.33) (0.049,191.33) (0.059,191.33) (0.069,191.33) (0.079, 191.33) (0.089,191.33) (0.099,191.33) (0.109,191.33) (0.119,191.33) (0.129,191.33) (0.139,191.33) };
\addlegendentry{\hspace{0.3cm}\textit{\textit{TGR} generation under \textit{CUC}}\hspace{1.5cm}\quad \hspace{1.6cm}}
\end{axis}
\end{tikzpicture}
\vspace{-0.4cm}
\caption{Influence of carbon tax rate on the \textit{EUC} and \textit{CUC} approaches}
\vspace{-0.4cm}
\label{sen1}
\end{figure}
While the simulation depicted in Fig. \ref{res} is performed for the specific carbon tax rate $\psi=\mathit{\$}0.07/kgCO_2e$, we also would like to examine the influence of the carbon tax rate on the quantifiable benefits of the \textit{EUC} approach. To this end, we study the sensitivity of the total costs to changes in the carbon tax rate. We vary the carbon tax rate $\psi$ from $\mathit{\$}0.009/kgCO_2e$ to the highest global carbon tax rate in 2018 of $\mathit{\$}0.139/kgCO_2e$ implemented in Sweden, in $\mathit{\$}0.010/kgCO_2e$ increments.\par
Fig. \ref{sen1} illustrates the energy generation from the \textit{TGR} $\gamma_1$ over the study period under the \textit{EUC} and \textit{CUC} approaches, as a function of $\psi$. Under \textit{EUC} approach, while the lowest simulated carbon tax rate resulted in a \textit{TGR} generation of 184.86 \textit{kWh}, the highest simulated carbon tax rate resulted in a significantly lower \textit{TGR} generation of 120 \textit{kWh}. The \textit{TGR} generation under the \textit{CUC} approach, however, does not vary with carbon tax rate and attains the constant value of 191.33 \textit{kWh}, because the \textit{CUC} approach does not consider the impact of carbon tax rate at the time of decision. The plots make clear that the carbon tax rate can effectively disincentivize \textit{TGR} generation under the \textit{EUC} approach in comparison with the \textit{CUC} solution.\par   
\enlargethispage{-0.1cm}
Fig. \ref{sen1} also presents the difference between the total costs obtained by the \textit{EUC} and \textit{CUC} approaches, \textit{i.e.}, the total costs under the \textit{CUC} approach minus the total costs under the \textit{EUC} approach, as a function of $\psi$. The results indicate that, for all considered carbon tax rates, the total costs under the \textit{EUC} approach are lower than those under the \textit{CUC} approach. We further observe that, as $\psi$ increases, the reduction in total costs under the \textit{EUC} approach vis-à-vis the \textit{CUC} approach also increases. This observation can be attributed to the fact that the \textit{EUC} approach optimizes the microgrid operation by taking into account the carbon tax payment based on varying carbon tax rates. The \textit{CUC} approach, however, can only evaluate the impact of increasing carbon tax rate ex-post, which inevitably results in higher costs with increasing carbon tax rates. These results underscore the importance of the explicit consideration of the monetary impacts of \textit{GHG} emissions by \textit{UC} approaches.
\section{Conclusion}\label{5}
In this paper, we propose a \textit{UC} approach that expressly assesses the \textit{GHG} emissions from \textit{TGR}s as well as their monetary impacts. The proposed \textit{EUC} approach enables the stipulation of a constraint on \textit{GHG} emissions from \textit{TGR}s over the study period and the ex-ante evaluation of carbon tax payment jointly with all other costs and benefits. The results indicate that the proposed approach yields lower costs than does the classical \textit{UC} approach. Further, it was observed that the \textit{TGR} operation could be attenuated via an economic mechanism, \textit{i.e.}, carbon tax rate, when the carbon tax payment is ex-ante evaluated. The performed sensitivity analysis provides valuable insights into the impact of carbon tax rate on the \textit{EUC} and \textit{CUC} approaches.\par
In our future studies, we plan to incorporate emissions trading schemes to the proposed \textit{EUC} approach. To evaluate a wider range of costs and benefits, we further plan to represent the participation of microgrids in wholesale energy and ancillary services markets under the \textit{EUC} approach.


\begin{thebibliography}{10}
\providecommand{\url}[1]{#1}
\csname url@samestyle\endcsname
\providecommand{\newblock}{\relax}
\providecommand{\bibinfo}[2]{#2}
\providecommand{\BIBentrySTDinterwordspacing}{\spaceskip=0pt\relax}
\providecommand{\BIBentryALTinterwordstretchfactor}{4}
\providecommand{\BIBentryALTinterwordspacing}{\spaceskip=\fontdimen2\font plus
\BIBentryALTinterwordstretchfactor\fontdimen3\font minus
  \fontdimen4\font\relax}
\providecommand{\BIBforeignlanguage}[2]{{%
\expandafter\ifx\csname l@#1\endcsname\relax
\typeout{** WARNING: IEEEtran.bst: No hyphenation pattern has been}%
\typeout{** loaded for the language `#1'. Using the pattern for}%
\typeout{** the default language instead.}%
\else
\language=\csname l@#1\endcsname
\fi
#2}}
\providecommand{\BIBdecl}{\relax}
\BIBdecl

\bibitem{eucbib:101}
D.~E. {Olivares} \emph{et~al.}, ``Trends in microgrid control,'' \emph{IEEE
  Transactions on Smart Grid}, vol.~5, no.~4, pp. 1905--1919, 2014.

\bibitem{eucbib:69}
J.~D. Glover, M.~S. Sarma, and T.~Overbye, \emph{Power System Analysis and
  Design}, 5th~ed.\hskip 1em plus 0.5em minus 0.4em\relax Stamford, CT: Cengage
  Learning, 2012, p.~32.

\bibitem{eucbib:63}
A.~J. Wood, B.~F. Wollenberg, and G.~B. Shebl{\'e}, \emph{Power Generation,
  Operation and Control}, 3rd~ed., Hoboken, NJ, 2014, pp. 148--153.

\bibitem{eucbib:79}
T.~A. {Nguyen} and M.~L. {Crow}, ``Stochastic optimization of renewable-based
  microgrid operation incorporating battery operating cost,'' \emph{IEEE
  Transactions on Power Systems}, vol.~31, no.~3, pp. 2289--2296, May 2016.

\bibitem{eucbib:62}
\BIBentryALTinterwordspacing
``State and trends of carbon pricing 2018,'' Tech. Rep., May 2018. [Online].
  Available: \url{https://openknowledge.worldbank.org/handle/10986/29687}
\BIBentrySTDinterwordspacing

\bibitem{eucbib:33}
X.~{Wu}, X.~{Wang}, and C.~{Qu}, ``A hierarchical framework for generation
  scheduling of microgrids,'' \emph{IEEE Transactions on Power Delivery},
  vol.~29, no.~6, pp. 2448--2457, Dec 2014.

\bibitem{eucbib:19}
A.~{Solanki}, A.~{Nasiri}, V.~{Bhavaraju}, Y.~L. {Familiant}, and Q.~{Fu}, ``A
  new framework for microgrid management: Virtual droop control,'' \emph{IEEE
  Transactions on Smart Grid}, vol.~7, no.~2, pp. 554--566, March 2016.

\bibitem{eucbib:95}
F.~{Li}, J.~{Qin}, and Y.~{Kang}, ``Closed-loop hierarchical operation for
  optimal unit commitment and dispatch in microgrids: A hybrid system
  approach,'' \emph{IEEE Transactions on Power Systems}, vol.~35, no.~1, pp.
  516--526, Jan 2020.

\bibitem{eucbib:96}
V.~{Sarfi}, I.~{Niazazari}, and H.~{Livani}, ``Multiobjective fireworks
  optimization framework for economic emission dispatch in microgrids,'' in
  \emph{2016 North American Power Symposium (NAPS)}, Sep. 2016, pp. 1--6.

\bibitem{eucbib:104}
Z.~{Ding} and W.~{Lee}, ``A stochastic microgrid operation scheme to balance
  between system reliability and greenhouse gas emission,'' in \emph{2015
  IEEE/IAS 51st Industrial Commercial Power Systems Technical Conference (I
  CPS)}, 2015, pp. 1--9.

\bibitem{eucbib:98}
\BIBentryALTinterwordspacing
O.~Yurdakul, ``Analysis and performance evaluation of coordinated transaction
  scheduling,'' M.S. thesis, University of Illinois at Urbana-Champaign,
  Urbana, IL, 2018. [Online]. Available:
  \url{http://hdl.handle.net/2142/101097}
\BIBentrySTDinterwordspacing

\bibitem{eucbib:83}
\BIBentryALTinterwordspacing
S.~{Greene} and A.~{Lewis}, ``{Global Logistics Emissions Council Framework for
  Logistics Emissions Accounting and Reporting},'' Tech. Rep., 2019. [Online].
  Available:
  \url{https://www.inlandwaterwaytransport.eu/wp-content/uploads/2019GLECFramework_Aug2019.pdf}
\BIBentrySTDinterwordspacing

\bibitem{eucbib:80}
\BIBentryALTinterwordspacing
{Pecan Street Inc.}, ``{Dataport},'' December 2019. [Online]. Available:
  \url{https://dataport.pecanstreet.org}
\BIBentrySTDinterwordspacing

\bibitem{eucbib:89}
\BIBentryALTinterwordspacing
{Consolidated Edison Company of New York, Inc}. (2019) {Time-of-Use Rates}.
  [Online]. Available:
  \url{https://www.coned.com/en/save-money/energy-saving-programs/time-of-use}
\BIBentrySTDinterwordspacing

\bibitem{eucbib:88}
D.~L. {Marino}, K.~{Amarasinghe}, and M.~{Manic}, ``Building energy load
  forecasting using deep neural networks,'' in \emph{IECON 2016 - 42nd Annual
  Conference of the IEEE Industrial Electronics Society}, Oct 2016, pp.
  7046--7051.

\end{thebibliography}
\end{document}